\documentclass[12pt, reqno, a4paper]{amsart}

\usepackage{amssymb}
\usepackage{color}
\definecolor{darkgreen}{rgb}{0,0.45,0} 
\usepackage[
colorlinks,citecolor=blue,linkcolor=blue, breaklinks]{hyperref}
\usepackage[matrix, arrow, curve]{xy}
\usepackage{bbm}
\usepackage{dsfont}
\usepackage{xcolor}

\usepackage[english]{babel}

\binoppenalty=\maxdimen
\relpenalty=\maxdimen

\theoremstyle{plain}
\newtheorem{theorem}{Theorem}[section]

\theoremstyle{remark}
\newtheorem{remark}[theorem]{Remark}

\theoremstyle{definition}
\newtheorem{example}[theorem]{Example}

\newtheorem{definition}[theorem]{Definition}

\numberwithin{equation}{section}

\DeclareMathOperator{\Id}{Id}

\DeclareMathOperator{\id}{id}
\DeclareMathOperator{\ch}{char}
\DeclareMathOperator{\Ker}{Ker}

\DeclareMathOperator{\sign}{sign}

\DeclareMathOperator{\Hom}{Hom}

\DeclareRobustCommand{\No}{\ifmmode{\nfss@text{\textnumero}}\else\textnumero\fi} 

\newbox\skewpullbackbox
\setbox\skewpullbackbox=\hbox{\xy 0;<1mm,0mm>: \POS(4,0)\ar@{-} (-4,0) \ar@{-} (8,4)
	\endxy}

\newbox\skwepullbackbox
\setbox\skwepullbackbox=\hbox{\xy 0;<1mm,0mm>: \POS(16,0)\ar@{-} (10,0) \ar@{-} (12,4)
	\endxy}

\newbox\ksewpullbackbox
\setbox\ksewpullbackbox=\hbox{\xy 0;<1mm,0mm>: \POS(0,-8)\ar@{-} (0,-4) \ar@{-} (4,-4)
	\endxy}

\newbox\pullbackbox
\setbox\pullbackbox=\hbox{\xy 0;<1mm,0mm>: \POS(4,0)\ar@{-} (0,0) \ar@{-} (4,4)
	\endxy}

\newbox\pullbackabox
\setbox\pullbackabox=\hbox{\xy 0;<1mm,0mm>: \POS(-4,-6)\ar@{-} (-8,-6) \ar@{-} (-4,-2)
	\endxy}

\newbox\pullbackbbox
\setbox\pullbackbbox=\hbox{\xy 0;<1mm,0mm>: \POS(-4,-4)\ar@{-} (-8,-4) \ar@{-} (-4,0)
	\endxy}

\newbox\pullbackcbox
\setbox\pullbackcbox=\hbox{\xy 0;<1mm,0mm>: \POS(4,-7)\ar@{-} (0,-7) \ar@{-} (4,-3)
	\endxy}

\newbox\pullbackdbox
\setbox\pullbackdbox=\hbox{\xy 0;<1mm,0mm>:\POS(-10,-6)\ar@{-} (-14,-6) \ar@{-} (-10,-2)
	\endxy}

\newbox\pushoutbox
\setbox\pushoutbox=\hbox{\xy 0;<1mm,0mm>: \POS(0,4)\ar@{-} (0,0) \ar@{-} (4,4)
	\endxy}

\newbox\pushoutabox
\setbox\pushoutabox=\hbox{\xy 0;<1mm,0mm>: \POS(2,8)\ar@{-} (2,4) \ar@{-} (8,8)
	\endxy}

\makeatletter
\@namedef{subjclassname@2020}{\textup{2020} Mathematics Subject Classification}
\makeatother

\begin{document}

\title[Hopf algebra polynomial identities in the Sweedler algebra]{Hopf algebra polynomial identities in the Sweedler algebra}

\author{Alexey Gordienko}

\email{alexey.gordienko@math.msu.ru}

\author{Maria Morozova}

\email{mariia.s.morozova@math.msu.ru}

\address{Department of Higher Algebra,
	Faculty of Mechanics and  Mathematics,
	M.\,V.~Lomonosov Moscow State University,
	Leninskiye Gory, d.1,  Main Building, GSP-1, 119991 Moskva, Russia }

\keywords{Polynomial identity, codimension, PROP, symmetric monoidal category, Hopf algebra, Sweedler algebra.}

\begin{abstract} Using Mac Lane's PROPs, one can introduce multilinear polynomial identities that involve all the operations of a Hopf algebra: the (co)multiplication, the (co)unit and the antipode. We find bases for polynomial identities and calculate their codimensions for the Sweedler algebra, treating it first as an ordinary associative algebra over a field and second as a Hopf algebra.	
\end{abstract}

\subjclass[2020]{Primary 16R10; Secondary 08B20, 16R50, 16T05, 18C10, 18M05.}

\thanks{The first author was supported by the Russian Science Foundation grant \No\,26-11-00001}

\maketitle



\section{Introduction}
Studying polynomial identities in an algebraic structure is an important aspect of understanding its properties.
Furthermore, when a basis for polynomial identities in a given algebra $A$ over a field is being calculated, certain numerical characteristics of polynomial identities, called codimensions, come into play naturally. The $n$th codimension $c_n(A)$ is precisely the dimension of the vector space of all linear maps $A^{\otimes n} \to A$ that can be realized by multilinear polynomials.
It is known~\cite[Section~6.2]{GiaZai} that there exists a tight relationship between the structure of~$A$ and the asymptotic behavior of  $c_n(A)$.

At the same time, \textit{$\Omega$-algebras} over a field~--- i.e., vector spaces $A$ with a set $\Omega$ of linear operations $A^{\otimes m} \to A^{\otimes n}$ where $m,n \in \mathbb Z_+$ (e.g., coalgebras, bialgebras, Hopf algebras, braided vector spaces, Yetter--Drinfel'd modules and other structures)~--- find their applications in many areas of mathematics and physics. 
 The traditional approach to defining a polynomial identity via elements of the corresponding free structure usually fails here because, for example, free coalgebras do not exist (only cofree coalgebras do, see~\cite[Section~6.4]{SweedlerBook}).
 One of the ways to overcome this difficulty is by treating a polynomial as a morphism in a suitable category called a theory.
 In~\cite{Lawvere} W.~Lawvere introduced \textit{Lawvere theories}, i.e., categories with the set of objects~$\mathbb{Z}_+$ and
 the products given by addition, and noticed that every variety of traditional algebraic structures (i.e., where the coarity $n$
 of all operations equals $1$) is equivalent to the category of product preserving functors from a suitable Lawvere theory to the category of sets. In order to include new algebraic structures, one has to replace a Lawvere theory with Mac Lane's PROP~\cite{MacLane65} (which is a symmetric strict monoidal category) or its braided analog as was done in~\cite{ASGordienko23} where  
 multilinear polynomial identities for coalgebras, bialgebras, Hopf algebras, braided vector spaces, Yetter--Drinfel'd modules and other structures were introduced. This approach makes it possible to define non-multilinear polynomial identities for traditional algebraic structures as well if one endows $A$ with two additional operations $\Delta$ and $\varepsilon$ defining the unique comonoid structure on $A$ and treats $A$ as an object in the category of sets.
 However, if $A$ is considered to be an object in the category of vector spaces, which is necessary when the coarity is allowed to be greater than $1$, one gets only multilinear polynomial identities. Moreover, if one treats an ordinary algebra $A$ over a field as a vector space, then one obtains not only ordinary multilinear polynomials, which take values in $A$, but also multilinear tensor polynomials, which take values in $A^{\otimes n}$. 
 
 Tensor polynomial identities  that are not necessarily multilinear for associative algebras were introduced in 2022 by F.~Huber and C.~Procesi~\cite{HuberProcesi}. Some results of the traditional PI-theory were generalized to tensor polynomial identities
 by D.\,Yu.~Novochadov in~\cite{Novochadov}. In addition, Novochadov~\cite[Proposition~1.4]{Novochadov} proved that if an associative algebra satisfies a nontrivial tensor polynomial identity, then it satisfies an ordinary polynomial identity.
 Moreover, in~\cite[Corollary~7.4]{Novochadov} he proved the analog of Amitsur's conjecture
 for tensor codimensions $c_{m,n}(A,\lbrace \mu \rbrace)$ of an arbitrary finite-dimensional associative algebra $A$ for fixed $n\in\mathbb N$
 and $m\to\infty$.
 
   In 2000, M.\,V.~Kochetov introduced~\cite{KochetovCoalgId} the notion of a (not necessarily multilinear) polynomial identity for a coalgebra $C$ as a polynomial identity for the algebra $C^*$ dual to~$C$. 
In 1996, M.~Markl used linear PROPs (which he simply called ``theories'') to study deformations of arbitrary $\Omega$-magmas~\cite{MarklPROPDef}. In 2002, T.~Pirashvili described explicitly the PROP corresponding to the notion of a bialgebra over a field~\cite{PirashviliPROPBialgebra}.
Polynomial identities in $\Omega$-algebras over a field, where the coarity of all operations equals $1$, were considered by Yu.\,A.~Bahturin, D.~Martinez Correa, D.~Diniz
and F.~Yasumura in~\cite{BahtYas1,BahtYas2,BahtMCDYas}.

So far, very few examples are known where a basis for polynomial identities has been found and their codimensions have been calculated explicitly. For instance, in~\cite{ASGordienko23}, this was done for the group Hopf algebra of the cyclic group of order $2$, a two-dimensional Yetter--Drinfel'd module and finite-dimensional vector spaces with no operations except swaps.
In Section~\ref{SectionSweedlerAsAHopfAlgebra} we find a basis for polynomial identities and compute their codimensions for the Sweedler algebra viewed as a Hopf algebra.  In Section~\ref{SectionSweedlerAsAnOrdinaryAlgebra}
we show that the ordinary polynomial identities and codimensions of the Sweedler algebra coincide with those
of the algebra $\mathrm{UT}_2(\mathbbm k)$ of upper triangular $2\times 2$ matrices,
while their tensor polynomial identities are different.

\section{Preliminaries}

\subsection{Ordinary polynomial identities in algebras over a field}\label{SubsectionOrdinaryPI}
We first recall the classical definition of polynomial identities and their codimensions in associative algebras~\cite{DrenKurs,GiaZai}.

Let $\mathbbm{k}$ be a field. Denote by $\mathbbm k \langle X \rangle$ the free non-unital associative algebra on
the countable set $X=\lbrace x_1, x_2, \ldots \rbrace$, i.e., the algebra of polynomials without a constant term with coefficients from $\mathbbm k$ in the non-commuting variables from $X$. Let $A$ be an associative $\mathbbm{k}$-algebra.
We say that $f\in \mathbbm k \langle X \rangle$ is an \textit{ordinary polynomial identity} in $A$ and write $f\equiv 0$
if $f(a_1, \ldots, a_n)=0$ for all $a_i \in A$ where $n\in\mathbb N$ is the number of variables that appear in $f$.
In other words, $f\equiv 0$ if and only if $\varphi(f)=0$ for every algebra homomorphism $\varphi \colon \mathbbm k \langle X \rangle \to A$. The set $\Id(A)$ of ordinary polynomial identities forms an ideal in $\mathbbm k \langle X \rangle$ that is invariant under all endomorphisms of $\mathbbm k \langle X \rangle$.

Consider the vector space $$P_n = \left\lbrace  \left.\sum\limits_{\sigma \in S_n} \alpha_\sigma\, 
x_{\sigma(1)} x_{\sigma(2)} \ldots x_{\sigma(n)} \right| \alpha_\sigma \in \mathbbm k \right\rbrace \subset \mathbbm k \langle X \rangle$$ of multilinear polynomials in $x_1, \ldots, x_n$.
(Here $S_n$ is the $n$th symmetric group.)

The number $c_n(A) := \dim \frac{P_n}{P_n \cap \Id(A)}$, $n\in\mathbb N$, is called the $n$th \textit{codimension}
of ordinary polynomial identities of $A$. Let  $\Hom_{\mathbbm k}(A^{{}\otimes n},A)$ be the space of all linear maps $A^{{}\otimes n} \to A$. Then $c_n(A)$ is precisely the dimension of the subspace of $\Hom_{\mathbbm k}(A^{{}\otimes n},A)$ consisting of all the maps that can be realized using the multiplication in $A$. A detailed study of the asymptotic behavior of the codimension sequence can be found in~\cite{GiaZai}.

\subsection{$\Omega$-magmas} Now we briefly recall the basic definitions from~\cite{ASGordienko23}. The main difference from~\cite{ASGordienko23}
is that, where applicable, we assume the base category to be symmetric rather than braided. 
For an account of monoidal categories,
we refer the reader to,  e.g.,~\cite{EGNObook}. The swap in a symmetric monoidal category is denoted by $\tau_{X,Y} \colon X \otimes Y \mathrel{\widetilde\to} Y \otimes X$.

Let $\Omega$ be a set together with maps $s,t \colon \Omega \to \mathbb Z_+$. We will refer to $\Omega$ as the \textit{signature}. 

\begin{definition} An \textit{$\Omega$-magma} in a monoidal category $\mathcal C$ is 
	an object $A$ endowed with morphisms $\omega_A \colon A^{\otimes s(\omega)} \to
	A^{\otimes t(\omega)}$ for every $\omega \in \Omega$. Here we use the convention that $A^{\otimes 0} := \mathbbm{1}$, the monoidal unit in $\mathcal C$.
\end{definition}

\begin{example}\label{ExampleMagma}
	Every magma (i.e., a set with a binary operation) is the same as an $\Omega$-magma in the category $\mathbf{Sets}$ of sets where $\Omega=\lbrace\mu \rbrace$ with $s(\mu)=2$ and $t(\mu)=1$.
\end{example}

\begin{example}\label{ExampleAlgebraOmega}
	Every algebra (not necessarily associative or unital)  over a field $\mathbbm{k}$
	is the same as an $\Omega$-magma in the category $\mathbf{Vect}_\mathbbm{k}$ of $\mathbbm k$-vector spaces where $\Omega=\lbrace\mu \rbrace$ with $s(\mu)=2$ and $t(\mu)=1$.
\end{example}

\begin{example}\label{ExampleHopfAlgebraOmega}
	Every Hopf algebra $H$ over a field $\mathbbm k$ is an $\Omega_\mathbf{Hopf}$-magma in $\mathbf{Vect}_{\mathbbm k}$
	where $\Omega_\mathbf{Hopf}=\lbrace \mu, u, \Delta, \varepsilon, S \rbrace$ are symbols for
	the multiplication $\mu_H \colon H \otimes H \to H$, the unit $u_H \colon \mathbbm k \to H$,
	the comultiplication $\Delta_H \colon H  \to H \otimes H$,
	the counit $\varepsilon_H \colon H  \to \mathbbm k$ and the antipode $S_H \colon H \to H$.
	We refer the reader to~\cite{DNR, Montgomery}
	for an account on Hopf algebras.
\end{example}

Let $A$ and $B$ be $\Omega$-magmas in a monoidal category $\mathcal C$. A morphism $f\colon A \to B$ is called an \textit{$\Omega$-magma homomorphism}
if for every $\omega \in \Omega$ the diagram below is commutative:
$$
\xymatrix{
	A^{\otimes s(\omega)} \ar[rr]^-{f^{\otimes s(\omega)}} \ar[d]_{\omega_A} && B^{\otimes s(\omega)} \ar[d]^{\omega_B}\\
	A^{\otimes t(\omega)} \ar[rr]^-{f^{\otimes t(\omega)}}  && B^{\otimes t(\omega)} 
}
$$

Denote by $\Omega\text{-}\mathsf{Magma}(\mathcal C)$ the category of $\Omega$-magmas in $\mathcal C$.

\subsection{$\Omega$-monomials}
Fix a set $\Omega$ and maps $s,t \colon \Omega \to \mathbb Z_+$. We now introduce expressions in the signature $\Omega$. Formally, they are morphisms in the symmetric strict monoidal category~$\mathcal M(\Omega)$ defined below.
The objects in~$\mathcal M(\Omega)$ are non-negative integers symbolizing monoidal powers of a single object.

Consider the free monoid $(M, \bullet)$ with the set of free generators $$\Omega \sqcup \lbrace \id_m, \tau_{m,n} \mid m,n \in \mathbb Z_+ \rbrace$$ and
define monoid homomorphisms $s \colon M \to\mathbb  Z_+$
and $t \colon M \to \mathbb Z_+$
extending maps $s, t \colon \Omega \to \mathbb Z_+$, such that
\begin{enumerate}
\item $s(\id_m)=t(\id_m)=m$;
\item $s(\tau_{m,n})=t(\tau_{m,n})=m+n$.
\end{enumerate}

Now consider the directed graph $\Gamma$ with the set of vertices $\mathbb Z_+$ and the set of edges $M$
where an edge $w \in M$ goes from the vertex $s(w)$ to the vertex $t(w)$.

Recall that an equivalence relation $\sim$ on hom-sets $X (a,b)$ of a category $X$ is a \textit{congruence} if
$g_1 f_1 \sim g_2 f_2$ for all $f_1, f_2 \in X (a,b)$, $g_1, g_2 \in X (b,c)$
and objects $a,b,c$ such that $f_1 \sim f_2$ and $g_1 \sim g_2$.
The \textit{factor category} $X/{\sim}$ is the category with the same objects as in $X$ and the hom-sets $X (a,b)/{\sim}$.

Let $X$ be the category where the objects are non-negative integers and the morphisms are finite paths
in $\Gamma$.
If $m,n \in \mathbb Z_+$, define $m \otimes n :=  m+n$.
 If $w_1 \ldots w_s \in X(m,n)$ and  $w_1'\ldots w_t' \in X(m',n')$
  for some $m,n,m',n' \in \mathbb Z_+$ and $w_1, \ldots, w_s, w_1', \ldots, w_t' \in M$, then define 
$w_1 \ldots w_s \otimes w_1' \ldots w_t' \in X(m+m',n+n')$ 
as $$(w_1 \bullet  \id_{n'}) \dots (w_s \bullet  \id_{n'})
 (\id_m \bullet  w_1') \ldots (\id_m \bullet  w_t').$$

Denote by $\mathcal M(\Omega)$ the factor category $X/{\sim}$ where  $\sim$ is the minimal congruence making 
\begin{enumerate}
\item $\id_m$ the identity morphisms of objects $m$;
\item $X/{\sim}$ a symmetric strict monoidal category with the swap $\tau_{m,n}$,
the monoidal unit~$0$ and the monoidal product $\otimes$.
\end{enumerate}

Then $\mathcal M(\Omega)$ is a \textit{PROP}, i.e., a strict symmetric monoidal category with the set of objects~$\mathbb Z_+$ and
the monoidal product given by addition~\cite[Section~24]{MacLane65}. The category~$\mathcal M(\Omega)$ is called the \textit{category of $\Omega$-monomials}.
The morphisms in $\mathcal M(\Omega)$ are called \textit{$\Omega$-monomials}.
 Define $s(f):= m$ and $t(f):= n$ for $f \in \mathcal M(\Omega)(m,n)$. Note that $1$ is an $\Omega$-magma in $\mathcal M(\Omega)$.

 \subsection{$\Omega$-polynomial identities with coefficients in a field}
 
 Let $\mathbbm k$ be a field. Recall that a category $\mathcal C$ is \textit{$\mathbbm k$-linear} if it is enriched over
 the category $\mathbf{Vect}_\mathbbm{k}$ of $\mathbbm k$-vector spaces, i.e., if all hom-sets of $\mathcal C$ are $\mathbbm k$-vector spaces
 and the composition of morphisms is $\mathbbm k$-bilinear. For $\mathbbm k$-linear monoidal categories we require that the monoidal product is $\mathbbm k$-bilinear too. A functor $F \colon \mathcal C \to \mathcal D$ between $\mathbbm k$-linear categories $\mathcal C$ and $\mathcal D$
 is called \textit{$\mathbbm k$-linear} if $F$ is a $\mathbbm k$-linear map on hom-sets.
  
 For $\Omega$-magmas in $\mathbbm k$-linear symmetric monoidal categories one can introduce polynomial identities with coefficients in $\mathbbm k$.
 
 Consider the symmetric monoidal category
 $\mathcal P(\mathbbm k, \Omega)$ where the objects, the monoidal product on them and the swap are the same as in $\mathcal M(\Omega)$
 and for every $m,n\in\mathbb{Z}_+$ the hom-set
 $\mathcal P(\mathbbm k, \Omega)(m,n)$ is the $\mathbbm k$-vector space with the formal basis $\mathcal M(\Omega)(m,n)$.
 The compositions and the monoidal product on morphisms are extended from $\mathcal M(\Omega)$ to
 $\mathcal P(\mathbbm k, \Omega)$ by the $\mathbbm k$-linearity.
 
 Then  $\mathcal P(\mathbbm k, \Omega)$ is a $\mathbbm k$-linear symmetric strict monoidal category, which we call the \textit{category of $\Omega$-polynomials with coefficients in $\mathbbm k$}. Morphisms in $\mathcal P(\mathbbm k, \Omega)$ are called \textit{$\Omega$-polynomials with coefficients in $\mathbbm k$}. 
  
 Let $A$ be an arbitrary $\Omega$-magma in a $\mathbbm k$-linear symmetric monoidal category $\mathcal C$. Denote by~$\mathcal E_A$ a $\mathbbm k$-linear symmetric strong monoidal functor $\mathcal P(\mathbbm k, \Omega) \to \mathcal C$ mapping
 \begin{enumerate}
 	\item  $1$ to $A$;
 	\item $\omega$ to $\omega_A$ for every $\omega \in \Omega$.
 \end{enumerate}
  Note that the functor $\mathcal E_A$, which is called the \textit{evaluation functor}, is completely determined by the choice of the natural isomorphism $\mathcal E_A(m)\otimes \mathcal E_A(n) \mathrel{\widetilde\to} \mathcal E_A(m+n)$ and an isomorphism $\mathbbm{1} \mathrel{\widetilde\to} \mathcal E_A(0)$ where $\mathbbm{1}$ is the monoidal unit in $\mathcal C$.
  
 An $\Omega$-polynomial $f$ is a \textit{polynomial identity in $A$ with coefficients in $\mathbbm k$}  if $\mathcal E_A(f) = 0$. In this case we say that $A$ satisfies $f\equiv 0$. Polynomial identities in $A$ with coefficients in $\mathbbm k$ are morphisms in the semicategory\footnote{The definition of a semicategory is identical to that of a category, except that the existence of identity morphisms is not required.} $\Id(A,\mathbbm k)$ where the set of objects is $\mathbb Z_+$ and $\Id(A,\mathbbm k)(m,n) := \Ker \left(\mathcal E_A \bigr|_{\mathcal P(\mathbbm k, \Omega)(m,n)}\right)$
 for every $m,n\in\mathbb Z_+$.
 
 Let $V$ be a set of $\Omega$-polynomials with coefficients in $\mathbbm k$. The full subcategory $\mathop\mathrm{Var}(V)$ of $\Omega\text{-}\mathsf{Magma}(\mathcal C)$, consisting of all $\Omega$-magmas satisfying $f\equiv 0$ for all $f\in V$, is called a \textit{variety} of $\Omega$-magmas. 
 
 \begin{example}
 	Hopf algebras form a variety $\mathrm{Var}(V_\mathbf{Hopf})$ for $\Omega_\mathbf{Hopf}=\lbrace \mu, u, \Delta, \varepsilon, S \rbrace$ and
 	\begin{equation}\label{EqVHopf}\begin{split}
 	V_\mathbf{Hopf}= \lbrace\mu(\mu \otimes \id_1)-\mu(\id_1 \otimes \mu ),\quad
 	\mu (u\otimes \id_1) - \id_1,\quad
 	\mu (\id_1 \otimes u) - \id_1,\\
 	(\Delta \otimes \id_1)\Delta - (\id_1 \otimes \Delta)\Delta,\quad
 	(\varepsilon \otimes \id_1)\Delta - \id_1, \quad
 	(\id_1 \otimes \varepsilon)\Delta - \id_1, \\
 	\varepsilon u - \id_0, \quad
 	\varepsilon\mu - \varepsilon \otimes \varepsilon, \\
 	\Delta u - u \otimes u,\quad
 	\Delta\mu - (\mu \otimes \mu)(\id_1 \otimes \tau_{1,1} \otimes \id_1)(\Delta \otimes \Delta),\\
 	\mu(S\otimes \id_1)\Delta - u\varepsilon,\quad
 	\mu(\id_1 \otimes S)\Delta - u\varepsilon\rbrace.
 	\end{split}
 	\end{equation}
 \end{example}
 
 When studying polynomial identities in an $\Omega$-magma that already lies in some variety $\mathop\mathrm{Var}(V)$, it is natural to identify with each other those polynomial identities that follow from each other by virtue of $V$.
 
 Recall that a \textit{monoidal ideal} $I$ in a $\mathbbm k$-linear monoidal category $X$ is
 a system of $\mathbbm k$-subspaces $I(a,b) \subseteq X(a,b)$, where $a,b$ are objects in $X$, such that
 \begin{enumerate}
 	\item $gf \in I(a,c)$ for all $f\in X(a,b)$, $g\in X(b,c)$ where either $f\in I(a,b)$
 	or $g\in I(b,c)$;
 	\item $f \otimes g \in I(a \otimes c, b \otimes d)$ for all $f\in X(a,b)$, $g\in X(c,d)$ where either $f\in I(a,b)$
 	or $g\in I(c,d)$
 \end{enumerate}
 for all objects $a,b,c,d$ in $X$. Then the \textit{factor category} $X/I$ is the category with the same objects as in $X$
 and the hom-sets $X(a,b)/I(a,b)$. It is easy to see that $X/I$ is again a $\mathbbm k$-linear monoidal category and the canonical functor $X \to X/I$ is a $\mathbbm k$-linear strict monoidal functor. If $X$ is symmetric, then the category $X/I$ and the functor $X \to X/I$ are symmetric too.
 
 \begin{example}
 	For every $\Omega$-magma $A$ the hom-sets $\Id(A,\mathbbm k)(m,n)$ of $\Id(A,\mathbbm k)$ form a monoidal ideal in $\mathcal P(\mathbbm k, \Omega)$.
 \end{example}
  
 Given a set $V$ of $\Omega$-polynomials with coefficients in $\mathbbm k$, we consider the minimal monoidal ideal $I$
 in $\mathcal P(\mathbbm k, \Omega)$ such that $V \subseteq \bigcup\limits_{m,n\in\mathbb Z_+} I(m,n)$.
 The category $\mathcal P_V(\mathbbm k, \Omega) := \mathcal P(\mathbbm k, \Omega)/I$ is called the category of \textit{$V$-relative
 	$\Omega$-polynomials with coefficients in $\mathbbm k$}. Then an $\Omega$-magma $A$ belongs to $\mathop\mathrm{Var}(V)$
 if and only if there exists a $\mathbbm k$-linear symmetric strong monoidal functor $F \colon \mathcal P_V(\mathbbm k, \Omega) \to \mathcal C$
 such that $F(1)=A$ and $F(\omega)=\omega_A$ for all $\omega \in \Omega$. The category
 $\mathcal P_V(\mathbbm k, \Omega)$ is a $\mathbbm k$-linear symmetric strict monoidal category.
  
 If an $\Omega$-magma $A$ belongs to $\mathop\mathrm{Var}(V)$ for some set $V$ of $\Omega$-polynomials with coefficients in~$\mathbbm k$, it is natural to define additional polynomial identities of~$A$ using $V$-relative
 $\Omega$-polynomials with coefficients in $\mathbbm k$. The functor $\mathcal P_V(\mathbbm k, \Omega) \to \mathcal C$
 induced by $\mathcal E_A$ is denoted again by $\mathcal E_A$.
 
 \begin{definition} Let $V$ be a set of $\Omega$-polynomials with coefficients in $\mathbbm k$.
 	A polynomial identity $f\equiv 0$, where $f$ is an $\Omega$-polynomial, \textit{follows} from $V$
 	 if the image of $f$ in $\mathcal P_V(\mathbbm k, \Omega)$ is zero.
 	Let $A$ be an $\Omega$-magma.  We say that $V$ \textit{generates}
 	polynomial identities in $A$  with coefficients in $\mathbbm k$ (or that $V$ is a \textit{basis} for them) if $A$
 	belongs to $\mathop\mathrm{Var}(V)$
 	and the functor $\mathcal E_A \colon \mathcal P_V(\mathbbm k, \Omega) \to \mathcal C$ is faithful. In other words, $V$ generates polynomial identities of $A$ if every polynomial identity $f \equiv 0$ holding in $A$ follows from~$V$.
 \end{definition}	
 
 \begin{example}\label{ExampleAssocAlgebra}
 	Let $s(\mu)=2$, $t(\mu)=1$, $V_\mathrm{assoc} := \lbrace \mu(\mu \otimes \id_1) -  \mu( \id_1 \otimes \mu )\rbrace$ and let $\mathbbm k$ be a field.
 	In other words, $V_\mathrm{assoc}$ consists of the associator $(xy)z-x(yz)$. Then $\mathcal P_{V_\mathrm{assoc}}(\mathbbm k, \lbrace \mu \rbrace)(n,1)$, where $n\in\mathbb N$, coincides with the vector space $P_n$ of associative multilinear polynomials in the variables $x_1, \ldots, x_n$. Under this identification,
 	polynomial identities in an arbitrary associative algebra $A$ from
 	the set $\Id(A, \mathbbm k)(n,1)$ correspond to ordinary multilinear polynomial identities in $A$
 	from the set $P_n \cap \Id(A)$
 	 (see Section~\ref{SubsectionOrdinaryPI}). The elements of $\Id(A, \mathbbm k)(m,n)$ for arbitrary $m,n\in\mathbb Z_+$
 	 are called \textit{multilinear tensor polynomial identities} of $A$.
 \end{example}	

 Consider an $\Omega$-magma $A$ in a $\mathbbm{k}$-linear symmetric monoidal category $\mathcal C$
 where $\mathbbm{k}$ is a field.  Let $m,n \in \mathbb{Z}_+$.  Then
  $$c_{m,n}(A, \Omega) := \dim \mathcal E_A\bigl(\mathcal P(\mathbbm{k}, \Omega)(m, n) \bigr)
  = \dim \mathcal P(\mathbbm{k}, \Omega)(m, n)/\Id(A, \mathbbm{k})(m, n)
  $$ is called the \textit{$(m,n)$-codimension} of polynomial identities of~$A$.
  If $A$ belongs to some variety $\mathrm{Var}(V)$, then $\mathcal P(\mathbbm{k}, \Omega)(m, n)$
  can be replaced in the definition of~$c_{m,n}(A, \Omega)$ by $\mathcal P_V (\mathbbm{k}, \Omega)(m, n)$, since
  $c_{m,n}(A, \Omega)$ equals the dimension of the vector space of all morphisms $A^{{}\otimes m} \to A^{{}\otimes n}$
realizable by swaps and operations $\omega_A$ for $\omega \in \Omega$.

 \begin{example} If $A$ is an ordinary algebra over a field $\mathbbm{k}$, then $c_{n,1}(A, \lbrace \mu \rbrace)$ equals the $n$th ordinary codimension $c_n(A)$ of polynomial identities of $A$ for $n\in \mathbb N$ (see Section~\ref{SubsectionOrdinaryPI} and Example~\ref{ExampleAssocAlgebra}).
 \end{example}	

	Let $H$ be a Hopf algebra over a field $\mathbbm k$.
While maps belonging to 
$\mathcal E_H\bigl(\mathcal M(\Omega)(m, n) \bigr)$ can be depicted using graphical calculus (see,  e.g.,~\cite{KasselQuantumGroups}), their compositions can also be represented in the usual functional notation.
For example, $\mu(\mu \otimes \id_1)$ can be written as $(xy)z$ or $(x_1 x_2)x_3$ and 
$(\Delta \otimes \id_1)\Delta$ can be written as $x_{(1)(1)}\otimes x_{(1)(2)}\otimes x_{(2)}$.

\section{Sweedler algebra as an associative algebra over a field}\label{SectionSweedlerAsAnOrdinaryAlgebra}

Recall that the \textit{Sweedler algebra} $H_4$ over a field $\mathbbm k$ is the Hopf algebra that is generated as a unital algebra by elements $c$ and $v$ where $c^2=1$, $v^2=0$, $vc=-cv$, $H_4 = \langle 1, c, v, cv \rangle_{\mathbbm k}$, $Sc=c$, $Sv=-cv$, $\Delta c = c\otimes c$, $\Delta v = c\otimes v + v\otimes 1$,
$\varepsilon(c)=1$, $\varepsilon(v)=0$.

We begin by treating the Sweedler algebra as an ordinary algebra.
In Theorem~\ref{TheoremSweedlerAsAnOrdinaryAlgebra} below, we show that if $\ch \mathbbm k = 0$, then the Sweedler algebra has the same ordinary polynomial identities and codimensions as the algebra 
$\mathrm{UT}_2(\mathbbm k)$ of upper triangular $2\times 2$ matrices (see, e.g.,~\cite[Theorem 4.1.5]{GiaZai}).

\begin{theorem}\label{TheoremSweedlerAsAnOrdinaryAlgebra}
Ordinary polynomial identities of the Sweedler algebra~$H_4$
as an associative algebra over a field~$\mathbbm k$, where $\ch \mathbbm k = 0$, 
are generated by \begin{equation}\label{EqTheoremSweedlerAsAnOrdinaryAlgebra}[x,y][z,t]\equiv 0
\end{equation}
 where $[x,y]:=xy-yx$. Their ordinary codimensions are equal to
$$c_n(H_4)=c_{n,1}(H_4, \lbrace \mu \rbrace)=2^{n-1}(n-2)+2,\quad n \in\mathbb N.$$
\end{theorem}
\begin{proof} Throughout the proof, the term ``polynomial identity'' will mean an ordinary polynomial identity of $H_4$ as an associative algebra over a field.
	
	In order to verify the multilinear polynomial identity~\eqref{EqTheoremSweedlerAsAnOrdinaryAlgebra},
	it is sufficient to check it on the basis elements $1,c,v,cv$. Every commutator $[x,y]$ takes values
	in the Jacobson radical $\langle v,cv\rangle_{\mathbbm k}$, which has zero multiplication. Hence the value of $[x,y][z,t]$
	is always zero and~\eqref{EqTheoremSweedlerAsAnOrdinaryAlgebra} follows.
	
	 We then apply the techniques from~\cite[Theorem 4.1.5]{GiaZai}.	
	Note that $$[x,y]w[z,t]=[x,y][w,[z,t]]+[x,y][z,t]w.$$ Hence the polynomial
	identity \begin{equation}\label{EqTheoremSweedlerAsAnOrdinaryAlgebraConsequence}[x,y]w[z,t]\equiv 0\end{equation} is a consequence of~\eqref{EqTheoremSweedlerAsAnOrdinaryAlgebra}.
	
	Consider an arbitrary multilinear monomial $x_{i_1} x_{i_2}\ldots x_{i_n}$. If $i_1 < i_2 < \ldots < i_n$,
	then $x_{i_1} x_{i_2}\ldots x_{i_n} = x_1 x_2 \ldots x_n$.
	If the sequence contains an inversion, we find $1\leqslant k \leqslant n-1$
	such that $i_1 < i_2 < \ldots < i_k$, $i_{k+1} < i_k$ and rewrite
	$$x_{i_1} x_{i_2}\ldots x_{i_n} =  x_{i_1} x_{i_2}\ldots x_{i_{k-1}}[x_{i_k}, x_{i_{k+1}}]
	x_{i_{k+2}} \ldots x_{i_n} + x_{i_1} x_{i_2}\ldots x_{i_{k-1}} x_{i_{k+1}} x_{i_k}
	x_{i_{k+2}} \ldots x_{i_n}.
	$$ 
	The second summand is treated analogously. The first summand
	is rewritten as \begin{equation*}\begin{split}x_{i_1} x_{i_2}\ldots x_{i_{k-1}}[x_{i_k}, x_{i_{k+1}}]
	x_{i_{k+2}} \ldots x_{i_n} = x_{i_1} x_{i_2}\ldots x_{i_{k-1}}[[x_{i_k}, x_{i_{k+1}}], 
	x_{i_{k+2}}] x_{i_{k+3}} \ldots x_{i_n} \\
	+ x_{i_1} x_{i_2}\ldots x_{i_{k-1}} x_{i_{k+2}}[x_{i_k}, x_{i_{k+1}}]
	x_{i_{k+3}} \ldots x_{i_n}.\end{split}\end{equation*}
	Using $yx=[y,x]+xy$, we continue the process above and move variables that stand outside the commutators from their right to their left. Ultimately, we express an arbitrary multilinear monomial as a linear combination of
	polynomials $x_1 x_2 \ldots x_n$ and $x_{i_1} x_{i_2}\ldots x_{i_k}[x_{j_1}, x_{j_2}, \ldots, x_{j_{n-k}}]$
	where $[y_1, \ldots, y_n] := [[\ldots[[y_1, y_2],y_3], \ldots], y_n]$.
	
	Using~\eqref{EqTheoremSweedlerAsAnOrdinaryAlgebra} and~\eqref{EqTheoremSweedlerAsAnOrdinaryAlgebraConsequence},
	we arrange the variables $x_{i_1} \ldots x_{i_k}$ in the ascending order.
	Moreover, by the Leibniz rule and the anticommutativity of the commutator we have $$[y,z,x]=[[y,z],x]=[[y,x],z]+[y,[z,x]]= -[x,y,z]+[x,z,y].$$
	Therefore we can always put a specific $x_s$ inside the inner commutator and make $j_1$ the least among $j_1, \ldots, j_{n-k}$.
	Finally, the polynomial identity $[[x,y],[z,t]]=[x,y][z,t]-[z,t][x,y]\equiv 0$,
	which is a consequence of~\eqref{EqTheoremSweedlerAsAnOrdinaryAlgebra} too,
	is used to make $j_3 < \ldots < j_{n-k}$.
	
	In other words, every multilinear monomial can be written modulo~\eqref{EqTheoremSweedlerAsAnOrdinaryAlgebra} as a linear combination
	of polynomials
	\begin{enumerate}
		\item $x_1 x_2 \ldots x_n$;
		\item $x_{i_1} x_{i_2}\ldots x_{i_k}[x_{j_1}, x_{j_2}, \ldots, x_{j_{n-k}}]$ where $$i_1 < i_2 < \ldots < i_k,\quad j_1 < j_2,\quad
		j_1 < j_3 < \ldots < j_{n-k}.$$
	\end{enumerate}

Suppose that a linear combination of the polynomials above is a polynomial identity in $H_4$
and not all coefficients are zero. Substituting $x_1 = x_2 = \ldots = x_n =1$,
we see that $x_1 x_2 \ldots x_n$ is the only polynomial that does not vanish. Therefore its coefficient is zero.

Now choose $x_{i_1} x_{i_2}\ldots x_{i_k}[x_{j_1}, x_{j_2}, \ldots, x_{j_{n-k}}]$ with a nonzero coefficient such that $k$ is maximal among the polynomials with nonzero coefficients.
Substitute \begin{equation}\label{EqSweedlerOrdinarySubstitution}
x_{i_1} = x_{i_2}=\ldots = x_{i_k} = 1,\quad
x_{j_2} = v,\quad x_{j_1}=x_{j_3}=\ldots = x_{j_{n-k}} = c.\end{equation} Then $x_{i_1} x_{i_2}\ldots x_{i_k}[x_{j_1}, x_{j_2}, \ldots, x_{j_{n-k}}]$ takes a nonzero value and all 
$$x_{\ell_1} x_{\ell_2}\ldots x_{\ell_s}[x_{m_1}, x_{m_2}, \ldots, x_{m_{n-s}}]$$ with $\lbrace i_1, i_2, \ldots, i_k \rbrace  \nsubseteq \lbrace \ell_1, \ell_2, \ldots, \ell_s \rbrace $ vanish
since $1$ occurs in the long commutator. 
The case $\lbrace i_1, i_2, \ldots, i_k \rbrace  \subsetneqq \lbrace \ell_1, \ell_2, \ldots, \ell_s \rbrace $
is excluded by the maximality of $k$.
If $$\lbrace \ell_1, \ell_2, \ldots, \ell_s \rbrace = \lbrace i_1, i_2, \ldots, i_k \rbrace,$$ then $j_1 = m_1$ since both $j_1$ and $m_1$ are the minimal numbers in $$\lbrace m_1, m_2, \ldots, m_{n-s} \rbrace = \lbrace j_1, j_2, \ldots, j_{n-k} \rbrace.$$

 The polynomial $x_{\ell_1} x_{\ell_2}\ldots x_{\ell_s}[x_{m_1}, x_{m_2}, \ldots, x_{m_{n-s}}]$ does not 
vanish only if $j_2 = m_2$ and $$x_{\ell_1} x_{\ell_2}\ldots x_{\ell_s}[x_{m_1}, x_{m_2}, \ldots, x_{m_{n-s}}] = x_{i_1} x_{i_2}\ldots x_{i_k}[x_{j_1}, x_{j_2}, \ldots, x_{j_{n-k}}]$$ since $v$ must appear in the inner commutator.
Hence $x_{i_1} x_{i_2}\ldots x_{i_k}[x_{j_1}, x_{j_2}, \ldots, x_{j_{n-k}}]$ is the only polynomial that does not vanish under
the substitution~\eqref{EqSweedlerOrdinarySubstitution}, its coefficient is zero and we get a contradiction. Therefore, the polynomials (1)--(2) are linearly independent modulo $\Id(H_4)$. Since every multilinear polynomial can be written modulo~\eqref{EqTheoremSweedlerAsAnOrdinaryAlgebra} as a linear combination
of polynomials (1)--(2), all multilinear polynomial identities in $H_4$ are consequences of~\eqref{EqTheoremSweedlerAsAnOrdinaryAlgebra}.
 Recall that over a field of characteristic $0$ all polynomial identities in an algebra are consequences of its multilinear polynomial identities~\cite[Theorem 1.3.8]{GiaZai}. Hence all polynomial identities of $H_4$ as an associative algebra over a field are generated by~\eqref{EqTheoremSweedlerAsAnOrdinaryAlgebra}.
 The codimension equals the number of polynomials (1)--(2). Hence
 \begin{equation*}\begin{split}c_n(A)= 1 + \sum\limits_{k=2}^n \tbinom{n}{k} (k-1)
 = 1 +  \sum\limits_{k=2}^n \tbinom{n}{k} k - \sum\limits_{k=2}^n \tbinom{n}{k} \\ = 1 +
 \sum\limits_{k=2}^n \frac{n!}{(n-k)!(k-1)!} - (2^n - n-1 ) =  1  + n \sum\limits_{\ell=1}^{n-1} \frac{(n-1)!}{(n-1-\ell)!\ell!}
 - 2^n+n+1 \\= 1 + n (2^{n-1}-1) - 2^n+n+1 = (n-2)2^{n-1} + 2,
 \end{split}\end{equation*}
 where we replaced $k$ with $(\ell+1)$.
\end{proof}	

Note that tensor polynomial identities in $H_4$ and $\mathrm{UT}_2(\mathbbm k)$ are different: 
$\mathrm{UT}_2(\mathbbm k)$, being three-dimensional, satisfies $\sum\limits_{\sigma\in S_4} (\sign \sigma) x_{\sigma(1)}\otimes x_{\sigma(2)} \otimes x_{\sigma(3)}\otimes x_{\sigma(4)}\equiv 0$, while $H_4$ does not since $\dim H_4 = 4$.

Since $\mathcal P_{V_\mathrm{assoc}}(\mathbbm k, \lbrace \mu \rbrace)(n,n)$, where $n\in\mathbb N$,
is the linear span of permutations and does not involve the multiplication $\mu$,
the codimensions $c_{n,n}(A,\lbrace \mu \rbrace)$ of every algebra $A$ over a field $\mathbbm k$ are completely determined by $d:=\dim A$. By~\cite[Theorem~5.1]{ASGordienko23},
$$ c_{n,n}(A,\lbrace \mu \rbrace)
\sim \alpha_d n^{-\frac{d^2-1}{2}}
d^{2n} \text{ as }n\to\infty$$
where $\alpha_d := \left(\frac{1}{\sqrt{2\pi}}\right)^{d-1}
\left(\frac{1}{2}\right)^{(d^2-1)/2} \cdot 1! \cdot 2! \cdot\ldots
\cdot (d-1)! \cdot d^{d^2/2}$.

By~\cite[Corollary~7.4]{Novochadov}, for every $n\in\mathbb N$
$$\lim\limits_{m\to\infty} \sqrt[m]{c_{m,n}(H_4,\lbrace \mu \rbrace)}
= \lim\limits_{m\to\infty} \sqrt[m]{c_{m,n}(\mathrm{UT}_2(\mathbbm k),\lbrace \mu \rbrace)} = 2.$$

Below are some tensor polynomial identities that hold both in $H_4$ and $\mathrm{UT}_2(\mathbbm k)$:
$$\sum\limits_{\sigma\in S_4} (\sign \sigma) x_{\sigma(1)} x_{\sigma(2)} \otimes x_{\sigma(3)}x_{\sigma(4)}\equiv 0,$$
$$[x,y,z]\otimes [t,u,w]\equiv [x,y,w]\otimes [t,u,z],$$
$$[x,y,z,t]\otimes [u,w]\equiv [x,y]\otimes [u,w,z,t].$$

\section{Sweedler algebra as a Hopf algebra}\label{SectionSweedlerAsAHopfAlgebra}

\subsection{Operators $e,p,q$ and polynomials $f_{u,w}$}\label{SubsectionOperatorsEPQandMonomialsFUW}

 Let $\mathbbm k$ be a field, $\ch \mathbbm k \ne 2$. We choose 
 two bases $\mathcal A = \lbrace 1, c, v, cv\rbrace$ and $\mathcal B =\lbrace 1, c-1, (c-1)v, v(c-1)\rbrace$
 of the Sweedler algebra~$H_4$ over~$\mathbbm k$. We treat each basis element as a single letter
 in the alphabet $\mathcal A$ or $\mathcal B$. In particular, we distinguish between
 the word $(c,v)$ over $\mathcal A$ of length $2$ and the word $cv$ of length $1$.  As usual, the notation $1^n$ means
 $(\underbrace{1,\ldots, 1}_n)$.
 Given a word $u \in \mathcal A^m$, $m\in \mathbb Z_+$, we denote by~$|u|_a$ the number of occurrences of $a\in \mathcal A$
 in~$u$. The same notation is used for words in~$\mathcal B$ as well.
 Finally, we introduce operators $e,p,q$, which will be used in the $V_{\mathbf{Hopf}}$-relative $\Omega_\mathbf{Hopf}$-polynomials below. Let $e x := x-\varepsilon(x)1$, $p := \frac{\id +S^2}{2}$, $q := \frac{\id -S^2}{2}$.
 
 \begin{remark}\label{RemarkSweedlerOperatorsEPQ}
 		 If we substitute $1$, $c$, $v$ or $cv$
 	for $x$ in $epx$, $qx$, $epx_{(1)}\otimes qx_{(2)}$ and $epx_{(2)}\otimes qx_{(1)}$,
 	then $epx$ is nonzero if and only if $x=c$, $qx$ is nonzero if and only if $x \in \lbrace v, cv\rbrace$,
 	$epx_{(1)}\otimes qx_{(2)}$ is nonzero if and only if $x=v$  and $epx_{(2)}\otimes qx_{(1)}$
 	is nonzero if and only if $x=cv$.
 \end{remark}
 
  Suppose $m,n \in \mathbb Z_+$. For every $u\in\mathcal A^m$ and $w\in\mathcal B^n$
  satisfying the conditions
  \begin{enumerate}
  	\item  $|u|_v + |u|_{cv} =  |w|_{(c-1)v} + |w|_{v(c-1)}$, 
  	\item $u=1^m$ if and only if $w=1^n$,
  \end{enumerate}
  we define polynomials $f_{u,w}\in \mathcal P_{V_{\mathbf{Hopf}}}(\mathbbm{k}, \Omega_\mathbf{Hopf})(m,n)$ below,
  which should satisfy the following properties:  \begin{enumerate}
  	\item  $f_{u, w}(\tilde u_1, \ldots, \tilde u_m)= 0$
  for all $\tilde u := (\tilde u_1, \tilde u_2,\ldots, \tilde u_m) \in \mathcal A^m$
  where $\tilde u \ne u$ and $|u|_1 \leqslant |\tilde u|_1$;
  \item  $f_{u, w}(u_1, \ldots, u_m)$ is a nonzero scalar multiple of $$w_1 \otimes w_2 \otimes\dots\otimes w_n$$
  where $(w_1, w_2,\ldots, w_n):=w$.
\end{enumerate}
  
  I. If $u=1^m$ and $w=1^n$,
  then $$f_{u,w}(x_1, \ldots, x_m):= \varepsilon(x_1)\varepsilon(x_2)\ldots \varepsilon(x_m) \underbrace{1\otimes \dots \otimes 1}_n.$$
  
  II. If $|u|_v = |u|_{cv}=0$, but $|u|_c > 0$, then 
\begin{equation*}\begin{split}f_{u,w}(x_1, \ldots, x_m) := \varepsilon(x_{i_1}) \ldots \varepsilon(x_{i_s})
\,\underbrace{1 \otimes \dots \otimes 1}_{k_1} \otimes (epx_{j_1(1)}) (epx_{j_2}) \ldots (epx_{j_t}) \\
\otimes \underbrace{1 \otimes \dots \otimes 1}_{k_2} \otimes\, epx_{j_1(2)} \otimes
\dots 
\otimes \underbrace{1 \otimes \dots \otimes 1}_{k_q} \otimes\, epx_{j_1(q)} \otimes \underbrace{1 \otimes \dots \otimes 1}_{k_{q+1}}
\end{split}
\end{equation*}
where $i_1< \ldots< i_s$ are the numbers of places in $u$ where $1$ occurs, $j_1 < \ldots < j_t$ are 
the numbers of places in $u$ where $c$ occurs and $w=:(1^{k_1},(c-1),1^{k_2},\ldots,(c-1),1^{k_q},(c-1),1^{k_{q+1}})$.

III. When $|u|_v + |u|_{cv} > 0$, the polynomial $f_{u,w}$ is defined in a slightly more complicated way.

  Let $u=(u_1, \ldots, u_m)$ and $w=(w_1, \ldots, w_n)$ where all $u_i\in\mathcal A$
and $w_j \in \mathcal B$. Denote $$k:= \min \bigl\lbrace 1\leqslant i \leqslant m \mathrel{\bigl|}
u_i \in \lbrace v,cv \rbrace \bigr\rbrace,\quad  
\ell:= \min \bigl\lbrace 1\leqslant j \leqslant n \mathrel{\bigl|}
w_j \in \lbrace (c-1)v, v(c-1) \rbrace \bigr\rbrace.$$

Define formal expressions $f_{u,w,j}$, $0\leqslant j \leqslant n$, which will be initial segments of $f_{u,w}$, inductively as follows.

Let $f_{u,w,0} := \varepsilon(x_{i_1}) \ldots \varepsilon(x_{i_s})$
where $\lbrace i_1, i_2, \ldots, i_s \rbrace := \lbrace 1\leqslant i \leqslant m \mid
u_i=1\rbrace$.

Suppose $f_{u,w,j}$ is already defined.

1. If $w_{j+1} = 1$, then $f_{u,w,j+1} := f_{u,w,j} \otimes 1$.

2. If $w_{j+1} = c-1$ and $u_k = v$, 
then $$f_{u,w,j+1} := f_{u,w,j} \otimes epx_{k(r)}$$ where $r:= 1+ |(w_1, \ldots, w_j)|_{c-1}$.

3. If $w_{j+1} = c-1$ and $u_k = cv$, 
then $$f_{u,w,j+1} := f_{u,w,j} \otimes epx_{k(r)}$$ where $r:= 3+ |(w_1, \ldots, w_j)|_{c-1}$.

4. Suppose $w_{j+1}\in \lbrace (c-1)v, v(c-1)\rbrace$ and $j+1 = \ell$.
Let $r:= |w|_{c-1}$.

a) If $w_\ell = (c-1)v$ and $u_k = v$, 
then $f_{u,w,\ell} := f_{u,w,\ell-1} \otimes (epx_{k(r+1)})\left(\prod\limits_{\substack{u_t = c,\\ 1\leqslant t \leqslant m}} epx_t\right) (qx_{k(r+2)})$.

b) If $w_\ell = (c-1)v$ and $u_k = cv$, 
then $f_{u,w,\ell} := f_{u,w,\ell-1} \otimes (epx_{k(2)})\left(\prod\limits_{\substack{u_t = c,\\ 1\leqslant t \leqslant m}} epx_t\right) (qx_{k(1)})$.

c) If $w_\ell = v(c-1)$ and $u_k = v$, 
then $f_{u,w,\ell} := f_{u,w,\ell-1} \otimes (qx_{k(r+2)})(epx_{k(r+1)})\left(\prod\limits_{\substack{u_t = c,\\ 1\leqslant t \leqslant m}} epx_t\right)$.

d) If $w_\ell = v(c-1)$ and $u_k = cv$, 
then $f_{u,w,\ell} := f_{u,w,\ell-1} \otimes (qx_{k(1)}) (epx_{k(2)})\left(\prod\limits_{\substack{u_t = c,\\ 1\leqslant t \leqslant m}} epx_t\right)$.

5. Suppose $w_{j+1}\in \lbrace (c-1)v, v(c-1)\rbrace$ and  $j+1 > \ell$. Let $$r:= |(w_1, \ldots, w_j)|_{(c-1)v} + 
|(w_1, \ldots, w_j)|_{v(c-1)}$$ and the number $1\leqslant t \leqslant m$ is defined by the condition that
$u_t \in \lbrace v, cv \rbrace$ and $$|(u_1, \ldots, u_{t-1})|_v + |(u_1, \ldots, u_{t-1})|_{cv} =  r.$$

a) If $w_{j+1} = (c-1)v$ and $u_t = v$, 
then $f_{u,w,j+1} := f_{u,w,j} \otimes (epx_{t(1)})(qx_{t(2)})$.

b) If $w_{j+1} = (c-1)v$ and $u_t = cv$, 
then $f_{u,w,j+1} := f_{u,w,j} \otimes (epx_{t(2)})(qx_{t(1)})$.

c) If $w_{j+1} = v(c-1)$ and $u_t = v$, 
then $f_{u,w,j+1} := f_{u,w,j} \otimes (qx_{t(2)})(epx_{t(1)})$.

d) If $w_{j+1} = v(c-1)$ and $u_t = cv$, 
then $f_{u,w,j+1} := f_{u,w,j} \otimes (qx_{t(1)})(epx_{t(2)})$.

Finally, let $f_{u,w}(x_1, \ldots, x_m) := f_{u,w,n}$.

\begin{example} For $u = (1,v,cv,c)$ and $w = \bigl(v(c-1), (c-1)v, c-1, 1\bigr)$  we have $k=2$ and $\ell = 1$. Hence
	$$f_{(1,v,cv,c),\bigl(v(c-1), (c-1)v, c-1, 1\bigr)}
	=  \varepsilon(x_1) (qx_{2(3)})(epx_{2(2)})(epx_4) \otimes (epx_{3(2)})(qx_{3(1)}) \otimes epx_{2(1)} \otimes 1.$$
\end{example}

\subsection{Calculation of codimensions}
Now we are able to provide a basis for polynomial identities in $H_4$ and calculate their codimensions. 

\begin{theorem} Let $\mathbbm k$ be a field, $\ch \mathbbm k \ne 2$. Then (multilinear) polynomial identities of $H_4$ as a Hopf algebra are generated by
\begin{equation}\label{EqSweedlerIdFirst}
-2qx \equiv (epx_{(1)})(qx_{(2)})+(epx_{(2)})(qx_{(1)})+(qx_{(2)})(epx_{(1)})+(qx_{(1)})(epx_{(2)}),
\end{equation}
\begin{equation}\label{EqSweedlerIdRearrangeFirst}epx\otimes epy \equiv epy \otimes epx,
\end{equation}
\begin{equation}(epx)(epy)\otimes epz \equiv epx \otimes (epy)(epz),
\end{equation}
\begin{equation}epx_{(1)}\otimes (epx_{(2)})(epy) \equiv -2 epx \otimes epy,
\end{equation}
\begin{equation} epx_{(1)} \otimes qx_{(2)} \otimes epy_{(1)} \otimes qy_{(2)}
\equiv epy_{(1)} \otimes qy_{(2)} \otimes epx_{(1)} \otimes qx_{(2)},
\end{equation}
\begin{equation} epx_{(2)} \otimes qx_{(1)} \otimes epy_{(2)} \otimes qy_{(1)}
\equiv epy_{(2)} \otimes qy_{(1)} \otimes epx_{(2)} \otimes qx_{(1)},
\end{equation}
\begin{equation} (epx_{(1)}) (qx_{(2)}) \otimes (epy_{(2)}) (qy_{(1)})
\equiv (epy_{(2)}) (qy_{(1)}) \otimes (epx_{(1)}) (qx_{(2)}),
\end{equation}
\begin{equation}\label{EqSweedlerIdRearrangeLast}  (qx_{(2)}) (epx_{(1)}) \otimes  (qy_{(1)}) (epy_{(2)})
\equiv  (qy_{(1)}) (epy_{(2)})\otimes  (qx_{(2)}) (epx_{(1)}),
\end{equation}
\begin{equation}\label{EqSweedlerIdZeroDoubleQFirst}qx_{(1)}\otimes qx_{(2)} \equiv 0,
\end{equation}
\begin{equation}(qx)(qy) \equiv 0,
\end{equation}
\begin{equation}\label{EqSweedlerIdZeroDoubleEPFromDifferentSidesOfQ}(epx)(qy)(epz) \equiv 0.
\end{equation}
\begin{equation}Spx \equiv px,
\end{equation}
\begin{equation}Sqx_{(1)} \otimes epx_{(2)} \equiv -(qx_{(1)})(px_{(2)}) \otimes epx_{(3)},
\end{equation}
\begin{equation}\label{EqSweedlerIdLast}Sqx_{(2)} \otimes epx_{(1)} \equiv (qx_{(3)})(px_{(2)}) \otimes epx_{(1)}.
\end{equation}
The images of polynomials $f_{u,w}$, defined in Section~\ref{SubsectionOperatorsEPQandMonomialsFUW}, form a basis of the $$\mathcal P_{\Id(H_4, \mathbbm{k})}(\mathbbm{k}, \Omega_\mathbf{Hopf})(m,n).$$
Moreover, $$c_{m,n}(H_4,\Omega_\mathbf{Hopf})= 1 + (2^m-1)(2^n-1) + 2^{m+n} \left(\tbinom{m+n}{n}-1\right)
\text{ for all }m, n\in\mathbb Z_+.$$
\end{theorem}
\begin{proof} Polynomial identities~\eqref{EqSweedlerIdFirst}--\eqref{EqSweedlerIdLast}
	are verified by substituting elements from the basis $\mathcal A = \lbrace 1, c, v, cv\rbrace$ and applying Remark~\ref{RemarkSweedlerOperatorsEPQ}.
	
	Note that the object $1$ in the category $\mathcal P_{V_\mathbf{Hopf}}(\mathbbm{k}, \Omega_\mathbf{Hopf})$ is a \textit{Hopf monoid}, i.e.,
	an object together with morphisms corresponding to the elements of the set $\Omega_\mathbf{Hopf}=\lbrace \mu, u, \Delta, \varepsilon, S \rbrace$ with the proper (co)arity,
	that satisfies the polynomial identities from $V_\mathbf{Hopf}$.
	By~\cite[Proposition~36]{Porst1}, the antipode in every Hopf monoid is a (co)monoid antihomomorphism.	
	Therefore, using the identities from $V_\mathbf{Hopf}$,
	every monomial from $\mathcal P_{V_\mathbf{Hopf}}(\mathbbm{k}, \Omega_\mathbf{Hopf})(m,n)$
	can be rewritten as a linear combination of monomials
	$$\varepsilon(x_{i_1})\ldots \varepsilon(x_{i_k}) \bigl(S^{\alpha_{11}}x_{\beta_{11}(\gamma_{11})}\bigr)
	\dots \bigl(S^{\alpha_{1\ell_1}}x_{\beta_{1\ell_1}(\gamma_{1\ell_1})}\bigr) \otimes \dots  \otimes
	\bigl(S^{\alpha_{s1}}x_{\beta_{s1}(\gamma_{s1})}\bigr)
	\dots \bigl(S^{\alpha_{s\ell_s}}x_{\beta_{s\ell_s}(\gamma_{s\ell_s})}\bigr),
	$$ where $\alpha_{ij} \in \mathbb Z_+$, $1\leqslant \beta_{ij} \leqslant m$
	and for every $1\leqslant i \leqslant m$ where $i \notin \lbrace i_1,  \ldots, i_k \rbrace$
	there exists $t\in \mathbb N$ such that $x_{i(j)}$ occurs in the monomial exactly once for every $1\leqslant j \leqslant t$
	and $x_{i(j)}$ does not occur for $j \geqslant t+1$.
	
	Now in every monomial we express $x_{i(j)} = \varepsilon(x_{i(j)}) 1 + epx_{i(j)} + qx_{i(j)}$ and
	apply~\eqref{EqSweedlerIdFirst} and~\eqref{EqSweedlerIdZeroDoubleQFirst}--\eqref{EqSweedlerIdLast}.
	Then in every resulting term all the antipodes $S$ disappear, i.e., become hidden in $p$ and $q$.
	By~\eqref{EqSweedlerIdZeroDoubleQFirst}--\eqref{EqSweedlerIdZeroDoubleEPFromDifferentSidesOfQ},
	in every tensor factor, expressions $qx$ occur at most once. Furthermore, for a given $1\leqslant i \leqslant m$,
	an expression $qx_{i(j)}$ appears at most once in the entire term too. Finally, in every tensor factor,
	expressions $epy$ cannot stand on both sides of $qx$, otherwise the term is zero.
	
	Applying~\eqref{EqSweedlerIdRearrangeFirst}--\eqref{EqSweedlerIdRearrangeLast}, we rearrange $epx$ and $qy$
	and get a linear combination of polynomials~$f_{u,w}$.
	 
	Now we prove that the images of $f_{u,w}$ in $\mathcal P_{\Id(H_4, \mathbbm{k})}(\mathbbm{k}, \Omega_\mathbf{Hopf})(m,n)$ are linearly independent.
	Suppose $f:=\sum\limits_{u,w} \alpha_{u,w} f_{u,w} (x_1, \ldots, x_m) \equiv 0$ in $H_4$
	for some $\alpha_{u,w} \in \mathbbm k$ not all equal to $0$.
	
	Choose $\alpha_{u,w}\ne 0$ with the maximal $|u|_1$ and substitute $x_1=u_1$, \ldots, $x_m=u_m$
	where $(u_1,u_2,\ldots, u_m):=u$. By Remark~\ref{RemarkSweedlerOperatorsEPQ}, all $f_{\tilde u, \tilde w}$ with $\alpha_{\tilde u,\tilde w}\ne 0$ vanish, except for $f_{u, \tilde w}$, which take nonzero values proportional to $\tilde w_1 \otimes \tilde w_2 \otimes\dots\otimes \tilde w_n$
	where $(\tilde w_1,\tilde w_2,\ldots, \tilde w_n):=\tilde w$. Note that $\tilde w_1 \otimes \tilde w_2 \otimes\dots\otimes \tilde w_n$
	are linearly independent for different $\tilde w$. Hence $\alpha_{u,w}= 0$, and we get a contradiction. Therefore,
	the images of $f_{u,w}$ indeed form a basis in $\mathcal P_{\Id(H_4, \mathbbm{k})}(\mathbbm{k}, \Omega_\mathbf{Hopf})(m,n)$.
	
	Since every $f \in \mathcal P_{V_\mathbf{Hopf}}(\mathbbm{k}, \Omega_\mathbf{Hopf})(m,n)$
	can be expressed as a linear combination of $f_{u,w}$ using only~\eqref{EqSweedlerIdFirst}--\eqref{EqSweedlerIdLast},
	the latter are indeed a basis for multilinear Hopf polynomial identities in~$H_4$. 
	
	Counting the number of the polynomials $f_{u,w}$ and using the Vandermonde identity $$ \sum\limits_{k=0}^{r} \tbinom{m}{k} \tbinom{n}{r-k}
	=\tbinom{m+n}{r},$$ we get
	\begin{equation*}\begin{split}c_{m,n}(H_4,\Omega_\mathbf{Hopf})=1 + (2^m-1)(2^n-1) + 2^{m+n} \sum\limits_{k=1}^{\min(m,n)} \tbinom{m}{k} \tbinom{n}{k}  \\
	=1 + (2^m-1)(2^n-1) + 2^{m+n} \sum\limits_{k=1}^{\min(m,n)} \tbinom{m}{k} \tbinom{n}{n-k}  \\
	 = 1 + (2^m-1)(2^n-1) + 2^{m+n} \left(\tbinom{m+n}{n}-1\right)
	\end{split}
	\end{equation*} for all $m, n\in\mathbb N$;
	$c_{0,n}(H_4,\Omega_\mathbf{Hopf})=c_{n,0}(H_4,\Omega_\mathbf{Hopf})=1$ for all $n\in \mathbb Z_+$.
\end{proof}	

\section*{AI usage}

The authors report generative AI was not used in their research or preparation of this manuscript.

\section*{Acknowledgments}

The first author is grateful to the Theoretical Physics and Mathematics Advancement Foundation ``BASIS''.
The study was conducted under the state assignment of Lomonosov Moscow State University.

\section*{Competing interests}

The authors confirm that there are no relevant financial or non-financial competing interests to report.


\begin{thebibliography}{99}

\bibitem{BahtYas1} Bahturin, Yu.\,A., Yasumura, F.  Distinguishing simple algebras by means of polynomial identities.
\textit{S\~ao Paulo J. Math. Sci.}, \textbf{13} (2019), 39--72.

\bibitem{BahtYas2} Bahturin, Yu.\,A., Yasumura, F. Graded polynomial identities as identities of universal algebras.
\textit{Lin. Algebra Appl.}, \textbf{562} (2019), 1--14.

\bibitem{BahtMCDYas} Bahturin, Yu.\,A., Martinez Correa, D., Diniz, D., Yasumura, F.  Polynomial identities of finite prime universal algebras. \texttt{arXiv:2512.06135 [math.RA] 5 Dec 2025}

\bibitem{DNR} D\u asc\u alescu, S., N\u ast\u asescu, C., Raianu, \c S.
Hopf algebras: an introduction. New York, Marcel Dekker, Inc., 2001.


\bibitem{DrenKurs} Drensky, V.\,S.
Free algebras and PI-algebras: graduate course in algebra.
Singapore, Springer-Verlag, 2000.

\bibitem{EGNObook} Etingof, P., Gelaki, S., Nikshych, D., Ostrik, V. Tensor categories.
\textit{AMS Mathematical Surveys and Monographs.} \textbf{205},
Providence, R.I., 2015, 362~pp.


\bibitem{GiaZai} Giambruno, A., Zaicev, M.\,V.
Polynomial identities and asymptotic methods.
\textit{AMS Mathematical Surveys and Monographs.} \textbf{122},
Providence, R.I., 2005, 352 pp.

\bibitem{ASGordienko23}
Gordienko, A.\,S.
On a general notion of a polynomial identity and codimensions.
\textit{J. Pure and Appl. Algebra}, \textbf{229}:1 (2025), 107814-1 -- 107814-27.

\bibitem{HuberProcesi} Huber, F., Procesi, C. Tensor polynomial identities. \textit{Israel J. Math.}, \textbf{247}:1 (2022), 125--147.

\bibitem{KasselQuantumGroups} Kassel, C. Quantum groups. \textit{Graduate Texts in Mathematics.}
\textbf{155}, Springer-Verlag, New York, 1995, 533 pp.

\bibitem{KochetovCoalgId} Kochetov, M.\,V. On identities for coalgebras and Hopf algebras.
\textit{Comm. Algebra}, \textbf{28}:3 (2000), 1211--1221.

\bibitem{Lawvere} Lawvere, F.\,W. Functorial semantics of algebraic theories. \textit{Proc. Nat. Acad.
Sci. U.S.A.}, \textbf{50} (1963), 869--873.

\bibitem{MacLane65} Mac Lane, S. Categorical algebra. \textit{Bulletin of the AMS}, \textbf{71}:1 (1965), 40--106.

\bibitem{MarklPROPDef} Markl, M. Cotangent cohomology of a category and deformations. 
\textit{J. of Pure and Applied Algebra}, \textbf{113} (1996), 195--218.


\bibitem{Montgomery} Montgomery, S. Hopf algebras and their actions on rings, CBMS Lecture Notes \textbf{82}, Amer. Math. Soc., Providence, RI, 1993.

\bibitem{Novochadov} Novochadov, D.\,Yu. Generalization of some results of the PI-theory to tensor identities. (2026). (Preprint, in Russian.)

\bibitem{PirashviliPROPBialgebra} Pirashvili, T. On the PROP corresponding to bialgebras. \textit{Cahiers de Topologie et G\'eom\'etrie Diff\'erentielle Cat\'egoriques}, \textbf{43}:3 (2002), 221--239.

\bibitem{Porst1}
Porst, H.-E. The formal theory of Hopf algebras. Part I: Hopf monoids in a monoidal category.  \textit{Quaestiones Mathematicae}, \textbf{38}:5 (2015), 631--682.

\bibitem{SweedlerBook} Sweedler M.E. Hopf Algebras, Benjamin, New York, 1969.

\end{thebibliography}
\end{document}